\documentclass{amsart}
\usepackage{amssymb, latexsym}
\usepackage{graphics}

\newcommand\C{{\mathbb C}} 
\newcommand\R{{{\mathbb R}}}

\newcommand\T{{{\mathbb T}}}
\newcommand\Schwartz{{{\mathcal S}}}
\newcommand\eps{\varepsilon}
\renewcommand\Re{{\operatorname{Re}}}

\newcommand\loc{{\operatorname{loc}}}
\newcommand\sym{{\operatorname{sym}}}

\theoremstyle{plain}
\newtheorem{theorem}{Theorem}
\newtheorem{remark}[theorem]{Remark}
\newtheorem{proposition}[theorem]{Proposition}
\newtheorem{lemma}[theorem]{Lemma}
\newtheorem{corollary}[theorem]{Corollary}

\numberwithin{equation}{section}
\numberwithin{theorem}{section}

\begin{document}

\title[Resonant decompositions and the $I$-method]{Resonant
  decompositions and the $I$-method for cubic nonlinear Schr\"odinger on $\R^2$}

\author{J. Colliander}
\thanks{J.C. is supported in part by N.S.E.R.C. Grant R.G.P.I.N. 250233-03.}
\author{M. Keel}
\thanks{M.K. was supported in part by N.S.F. Grant DMS-0303704
   and by the McKnight and Sloan Foundations.}
\author{G. Staffilani}
 \thanks{G.S. is supported in part by N.S.F. Grant DMS-0602678.}
 \author{H. Takaoka}\thanks{H.T is supported in part by J.S.P.S Grant No. 19740074.}
 \author{T. Tao}
 \thanks{T.T. is supported in part by a grant from the MacArthur
   Foundation. }
\date{20 April 2006}

\begin{abstract}
The initial value problem for the cubic defocusing nonlinear Schr\"odinger equation 
$i \partial_t u + \Delta u = |u|^2 u$ on the plane is shown to be globally well-posed for initial data in
$H^s ( \R^2)$ provided $s>1/2$. The proof relies upon an almost
conserved quantity constructed using multilinear correction
terms. The main new difficulty is to control the contribution of resonant interactions to these
correction terms.  The resonant interactions are 
significant due to the multidimensional setting of the problem and some orthogonality issues
which arise.
\end{abstract}

\maketitle

\section{Introduction}\label{section_1}

We consider the Cauchy problem for the cubic defocusing nonlinear Schr\"odinger (NLS) equation 
\begin{eqnarray}\label{eq:NLS}
\left\{
\begin{array}{l}
i\partial_t u+\Delta u=|u|^2u,\\
u(0,x)=u_0(x) \in H^s_x(\R^2),
\end{array}
\right.
\end{eqnarray}
in a Sobolev space $H^s_x(\R^2)$,
where the unknown function $u: J \times \R^2 \longmapsto \C$ is a strong solution to \eqref{eq:NLS} on a time interval $J \subset \R$ in the sense that $u \in C^0_{t,\loc} H^s_x(\R^d)$ and $u$ obeys the integral equation
$$ u(t) = e^{it\Delta} u_0 - i\int_0^t e^{i(t-t')\Delta}[|u|^2 u(t')]\ dt'$$
for $t \in J$.  Here of course the propagators $e^{it\Delta}$ are defined via the Fourier transform
$$ \hat f(\xi) := \int_{\R^2} e^{-i x \cdot \xi} f(x)\ dx$$
by the formula
$$ \widehat{e^{it\Delta} f}(\xi) := e^{-it|\xi|^2} \hat f(\xi)$$
and the Sobolev space $H^s_x(\R^2)$ is similarly defined via the Fourier transform using the norm
$$ \|f\|_{H^s_x(\R^2)} := \| \langle \xi \rangle^s \hat f(\xi) \|_{L^2_\xi(\R^2)}$$
where $\langle \xi \rangle := (1 + |\xi|^2)^{1/2}$.  For later use we shall also need the homogeneous Sobolev norms
$$ \|f\|_{\dot H^s_x(\R^2)} := \| |\xi|^s \hat f(\xi) \|_{L^2_\xi(\R^2)}.$$
We are interested primarily in the global-in-time problem, in which we allow $J$ to be the whole real line $\R$.

Both the local and global-in-time Cauchy problems for this NLS equation \eqref{eq:NLS} have attracted a 
substantial literature \cite{tsutsumi}, \cite{cw}, \cite{kato}, \cite{gv}
\cite{Staff}, \cite{Bo1}, \cite{Bo2}, \cite{CKSTT2}, 
\cite{planchon2}, \cite{cvv},  \cite{birnir}, \cite{cct}.  The local well-posedness theory is now well 
understood; in particular, one has local well-posedness in $H^s_x(\R^2)$ for all $s \geq 0$, 
and if $s$ is strictly positive then a solution can be continued unless the $H^s_x(\R^2)$ norm 
of the solution goes to infinity at the blowup time (see e.g. \cite{caz}, \cite{tao}).  Also, due to 
the smooth nature of the nonlinearity, any local $H^s(\R^2)$ solution can be expressed as the 
limit (in $C^0_{t,\loc} H^s_x$) of smooth solutions.  For $s < 0$ the solution map ceases to 
be uniformly continuous \cite{cct} and may possibly even be undefined, though it is known 
that well-posedness can be recovered for other spaces 
rougher than $L^2_x(\R^2)$ 
\cite{planchon2}, \cite{cvv}.  The 
space $L^2_x(\R^2)$ is the critical space for this equation, as it is invariant 
under the scaling symmetry
\begin{equation}\label{scaling}
u(t,x) \mapsto \frac{1}{\lambda} u(\frac{t}{\lambda^2}, \frac{x}{\lambda})
\end{equation}
of \eqref{eq:NLS}.

Now we turn attention to the global-in-time well-posedness problem.
Based on the local well-posedness theory, standard limiting arguments,
and the time reversal symmetry $u(t,x) \mapsto \overline{u(-t,x)}$, 
global well-posedness of \eqref{eq:NLS} for arbitrarily large 
data\footnote{Global well-posedness and even scattering is known when
  the mass $\| u_0 \|_{L^2_x(\R^2)}$ is sufficiently small (see
  e.g. \cite{caz}, \cite{tao}), or if suitable decay conditions
  (e.g. $x u_0 \in L^2_x(\R^2)$ are also imposed on the initial data
  \cite{tsutsumi}).  Our interest here however is in the large data
  case with no further decay conditions beyond the requirement that
  $u_0$ lie in $H^s_x(\R^2)$.} in $H^s_x(\R^2)$ for some $s > 0$
follows if an \emph{a priori bound} of the form
\begin{equation}\label{apriori}
\| u(T) \|_{H^s_x(\R^2)} \leq C( s, \|u_0\|_{H^s_x(\R^2)}, T )
\end{equation}
can be established for all times $0 < T < \infty$ and all smooth-in-time, Schwartz-in-space solutions $u: [0,T] \times \R^2 \to \C$, where the right-hand side is some finite quantity depending only upon $s$, $\|u_0\|_{H^s_x(\R^2)}$, and $T$.  Thus we shall henceforth restrict attention to such smooth solutions, which will in particular allow us to justify all formal computations, such as verification of conservation laws.

As is well known, the equation NLS enjoys two useful conservation laws, the \emph{energy conservation law}
\begin{eqnarray}\label{eq:energy}
E(u(t)):=\int_{\R^{2}}\frac{1}{2}|\nabla u(t,x)|^2+\frac{1}{4}|u(t,x)|^4\,dx = E(u_0).
\end{eqnarray}
and the \emph{mass conservation law}
\begin{equation}\label{masscons}
\| u(t) \|_{L^2_x(\R^2)} = \| u_0 \|_{L^2_x(\R^2)}.
\end{equation}
From these laws one easily establishes \eqref{apriori} for $s=1$ (with
bounds uniform in $T$), and with some additional arguments one can
then deduce the same claim for $s > 1$ (with the best known bounds
growing polynomially in $T$; see \cite{Staff}, \cite{CDKS}).  The mass conservation law \eqref{masscons} also gives \eqref{apriori} for $s=0$, but unfortunately this does not immediately imply any result for $s > 0$ except in the small mass 
case\footnote{In order to establish a global well-posedness result in $L^2_x(\R^2)$, it is 
instead necessary to obtain an a priori \emph{spacetime} bound such as 
$\| u \|_{L^4_{t,x}([0,T] \times \R^2)} \leq C( \|u_0\|_{L^2_x(\R^2)}
)$.  See \cite{keraani}, \cite{blue}, \cite{tao-conf} for further discussion.}.  

It is conjectured that the equation \eqref{eq:NLS} is globally well-posed in $H^s_x(\R^2)$ for all $s \geq 0$, and in particular \eqref{apriori} holds for all $s > 0$.  This conjecture remains open (though in the radial case, the higher dimensional analogue of this conjecture has recently been settled in \cite{tvz}).  However, there has been some progress in improving the $s \geq 1$ results mentioned earlier.  The first breakthrough was by Bourgain \cite{Bo1}, \cite{Bo2}, who established \eqref{apriori} (and hence global well-posedness in $H^s_x(\R^2)$) for all $s > 3/5$, using what is now referred to as the \emph{Fourier truncation method}.  

In \cite{CKSTT2} the bound \eqref{apriori} was established for all $s > 4/7$, using the ``$I$-method'' developed by the authors in \cite{ckstt}, \cite{CKSTT1} (see also \cite{keel:wavemap}).  
The main result of this paper is the following improvement:

\begin{theorem}[Main theorem]\label{main}  The bound \eqref{apriori} holds for all $s > 1/2$.  In particular, the Cauchy problem \eqref{eq:NLS} is globally well-posed in $H^s_x(\R^2)$ for all $s > 1/2$.
\end{theorem}

Our arguments refine our previous analysis in \cite{CKSTT2} by adding a ``correction term'' to
a certain modified energy functional $E(Iu)$, as in \cite{CKSTT1} or \cite{CKSTT3}, in order to damp out some oscillations in that functional; also, we establish some more refined estimates on the multilinear symbols appearing in those integrals.  
The main new difficulty is that, due to the non-integrability and multidimensional setting of this equation (in contrast\footnote{The equation considered in \cite{CKSTT3} was also non-integrable, but because it was one-dimensional there was still enough cancellation to prevent the contribution of the resonant interactions from becoming singular.} to \cite{CKSTT1}), 
the direct analogue of the correction terms used in \cite{CKSTT1}, \cite{CKSTT3} contains a singular symbol and 
is thus intractable to estimate.  We get around this new difficulty by truncating the correction term to 
non-resonant interactions, and dealing with the resonant interactions separately by
some advanced estimates of $X^{s,b}$ type.  This method seems quite
general and should lead to improvements in global well-posedness
results for other non-integrable evolution equations which are
currently obtained by the ``first-generation'' $I$-method
(i.e. without correction terms). A resonant decomposition similar to
that employed here appeared previously in the work  \cite{B:IMethod}, and more recently in \cite{BDGSz}.

Inserting the above theorem into the results of \cite{blue} (which employ the pseudo-conformal transform) we conclude that the 
equation \eqref{eq:NLS} is globally well-posed with scattering when the initial data obeys $\langle x \rangle^s u_0 \in L^2_x(\R^2)$ for any $s > 1/2$.

During the preparation of this manuscript, we learned that Fang and
Grillakis \cite{gf} had also obtained Theorem \ref{main}, in fact for
$s \geq 1/2$, by a different
method based upon a new type of Morawetz inequality. The
Fang-Grillakis interaction Morawetz estimate has recently \cite{CGTz} been
improved and combined with the $I$-method (following the general scheme from
\cite{CKSTT:CPAM}) to prove that \eqref{eq:NLS} is globally well-posed
in $H^s$ for $s > 2/5.$ The techniques leading to the improved energy 
increment control obtained in this paper (see \eqref{acl} which is
$N^{-1/2}$ better than what was obtained in \cite{CKSTT2} and used in
\cite{CKSTT:CPAM}, \cite{CGTz}) may also improve the ``almost
Morawetz'' increment in \cite{CGTz} by $N^{-1/2}$. Such an improvement
combined with \eqref{acl} would improve the global well-posedness
result to $s>4/13$. The arguments in \cite{gf}, \cite{CGTz} are based
on Morawetz inequalities and are thus restricted to the defocusing
case. Provided the mass of the initial data is less than the mass of
the ground state, Theorem \ref{main} also holds true for the focusing
analog of \eqref{eq:NLS} (see Remark \ref{focusingremark}
below). The focusing
problem is expected to be globally well-posed for $L^2$ initial data
with mass less than the ground state mass.


\subsection{Acknowledgements}  We thank Tristan Roy for detailed comments and corrections, and Manoussos Grillakis and Yung-Fu Fang for sharing their preliminary manuscript \cite{gf}.

\section{Setting up the $I$-method}

We now begin the proof of Theorem \ref{main}.  As in all other applications of the $I$-method, we will reduce matters to one of constructing a certain modified energy functional $\tilde E(u(t))$ and demonstrating that it has certain almost conservation properties.

By the discussion in the introduction, it suffices to prove \eqref{apriori} in the range $1/2 < s < 1$.  Henceforth we fix $s$.  We adopt the usual notation that $X \lesssim Y$ or $Y \gtrsim X$ denotes an estimate of the form $X \leq C(s) Y$, for some constant
$0 < C(s) < \infty$ depending only on $s$.  We also write $X \sim Y$ for $X \lesssim Y \lesssim X$, and $X = O(Y)$ for $|X| \lesssim Y$.

We will use exponents $a+$ and $a-$ to denote $a+\eps$ and $a-\eps$ for arbitrarily small exponents $\eps > 0$, and allow the implied constants in the $\lesssim$ notation to depend on $\eps$.  Thus for instance if we write $X \lesssim N^{1+} Y$, this means that for every $\eps$ there exists a constant $C(s,\eps)$ such that $X \leq C(s,\eps) N^{1+\eps} Y$.

Let $N \gg 1$ be a large parameter to be chosen later (it will eventually depend on $T$, $s$, and the size of the initial data $u_0$).
We define the Fourier multiplier $I = I_N$ by
$$ \widehat{Iu}(\xi) := m(\xi) \hat u(\xi)$$
where $m$ is a smooth non-negative radial symbol which equals $1$ when $|\xi| \leq N$, equals $(|\xi|/N)^{s-1}$ for $|\xi| \geq 2N$, and smoothly interpolates between the two in the region $N \leq |\xi| \leq 2N$.  We shall abuse notation and write $m(|\xi|)$ for $m(\xi)$, thus for instance $m(N)=1$.

The ``first-generation'' $I$-method revolves around the modified energy 
\begin{equation}\label{eiu}
E(Iu(t)) = \int_{\R^2} \frac{1}{2} |\nabla Iu(t,x)|^2 + \frac{1}{4} |Iu(t,x)|^4\ dx \sim
\| Iu \|_{\dot H^1_x(\R^2)}^2 + \| Iu \|_{L^4_x(\R^2)}^4,
\end{equation} 
and in particular establishing an almost conservation law for this quantity.  Here, we shall introduce a slight variant $\tilde E(u(t))$ of $E(Iu(t))$ and establish an almost conservation law for that quantity instead.  More precisely, we shall show the following:

\begin{theorem}[Existence of an almost conserved quantity]\label{eacq}  There exists a functional $\tilde E = \tilde E_N: \Schwartz_x(\R^2) \to \R$ defined on Schwartz functions $u \in \Schwartz_x(\R^2)$ with the following properties.
\begin{itemize}
\item (Fixed-time bounds) For any $u \in \Schwartz_x(\R^2)$, we have
\begin{equation}\label{fixed-bound}
 | E(Iu) - \tilde E(u) | \lesssim N^{-1+} \| Iu \|_{H^1_x(\R^2)}^4.
\end{equation}
\item (Almost conservation law)  If $\|u_0\|_{L^2_x(\R^2)} \leq A$ and $E(Iu_0) \leq 1$, and $u$ is a smooth-in-time, Schwartz-in-space solution to \eqref{eq:NLS} on a time interval $[0,t_0]$, then if $t_0$ is sufficiently small depending on $A$, we have
\begin{equation}\label{acl}
|\tilde E(u(t_0)) - \tilde E(u_0)| \lesssim C(A) N^{-2+}
\end{equation}
for some constant $C(A)$ depending only on $A$.
\end{itemize}
\end{theorem}

\begin{remark} The precise value of the exponent $-1+$ in \eqref{fixed-bound} is not particularly important; any negative exponent would have sufficed here.  However, the exponent $-2+$ in \eqref{acl} is directly tied to the restriction $s > 1/2$ in our main theorem.  More generally, an exponent of $-\alpha+$ in this almost conservation law translates to a constraint $s > 2/(2+\alpha)$.  In \cite{CKSTT2}, the first-generation modified energy $E(Iu)$ was shown to obey an almost conservation law with $\alpha = 3/2$, which ultimately led to the constraint $s > 4/7$.  Note that in order to get arbitrarily close to the scaling exponent $s=0$, one would need $\alpha$ to be arbitrarily large, which looks unlikely to be achieved with this method due to the lack of complete integrability.
\end{remark}

We shall prove Theorem \ref{eacq} in later sections.  For the remainder of this section, we show how Theorem \ref{eacq} implies Theorem \ref{main}.

\begin{proof}[Proof of Theorem \ref{main} assuming Theorem \ref{eacq}]  Fix $u$, $u_0$, $T$ as in Theorem \ref{main}, and write $A := 1 + \|u_0\|_{H^s_x(\R^d)}$.  We let $\lambda \geq 1$ be a scaling parameter to be chosen shortly, and define the rescaled solution $u^{(\lambda)}: [0, \lambda^2 T] \times \R^2 \to \C$ as per \eqref{scaling}, thus
$$ u^{(\lambda)}(t, x) := \frac{1}{\lambda} u(\frac{t}{\lambda^2}, \frac{x}{\lambda}).$$
Now let $N \gg 1$ also be a parameter to be chosen later (it will depend on $T$ and $A$).  A simple computation (see equation (3.10) of \cite{CKSTT2}) shows that
$$ E( I u^{(\lambda)}(0) ) \lesssim N^{2-2s} \lambda^{-2s} A^4.$$
Thus we can arrange
\begin{equation}\label{iuo}
 E( I u^{(\lambda)}(0) ) \leq 1/3
 \end{equation}
by choosing
\begin{equation}\label{lambda-def}
\lambda := C(s,A) N^{(1-s)/s}
\end{equation}
for a suitable quantity $C(s,A)$.  Also, from mass conservation (and scale-invariance) we also
know that
\begin{equation}\label{massive}
 \| Iu^{(\lambda)}(t) \|_{L^2_x(\R^2)} \leq A.
\end{equation}
From \eqref{fixed-bound}, \eqref{eiu} we conclude that
\begin{equation}\label{fixbo}
| E( I u^{(\lambda)}(t) ) - \tilde E( u^{(\lambda)}(t) ) | \lesssim N^{-1+} ( A^4 + E( I u^{(\lambda)}(t) )^2 ).
\end{equation}
We now claim that (for $\eps$ chosen suitably small, and for $N$ chosen suitably large) 
\begin{equation}\label{iulam}
 E( I u^{(\lambda)}(\lambda^2 T) ) \leq 2/3.
\end{equation}
To see this, suppose for contradiction that this were not the case; then there exists $0 < T' < \lambda^2 T$ such that
\begin{equation}\label{iueps}
 E( I u^{(\lambda)}(T') ) = 2/3
\end{equation}
but that
$$ E( I u^{(\lambda)}(t) ) \leq 2/3 \hbox{ for all } 0 \leq t \leq T'.$$
Applying \eqref{fixbo} we conclude (if $N$ is sufficiently large depending on $A$) that
$$ \tilde E( u^{(\lambda)}(t) ) \leq 1 \hbox{ for all } 0 \leq t \leq T'.$$
Applying \eqref{acl} repeatedly (and exploiting time translation invariance), we conclude that
$$ |\tilde E( u^{(\lambda)}(T') ) - \tilde E( u^{(\lambda)}(0) )| \lesssim C(A) N^{-2+} T' \leq
C(A) N^{-2+} \lambda^2 T$$
and hence by \eqref{fixbo}
$$  |E( Iu^{(\lambda)}(T') ) - E( Iu^{(\lambda)}(0) )| \lesssim C(A) N^{-2+} \lambda^2 T + N^{-1+} A^4.$$
From \eqref{lambda-def} and the hypothesis $s > 1/2$, we see that the net powers of $N$ on the right are negative.  Thus we can choose $N$ so large (depending on $A,T$) that
$$  |E( Iu^{(\lambda)}(T') ) - E( Iu^{(\lambda)}(0) )| < 1/3.$$
But this contradicts \eqref{iuo}, \eqref{iueps}.  Thus \eqref{iulam} must hold.  From this, \eqref{massive}, and some Fourier analysis we deduce
$$ \| u^{(\lambda)} \|_{\dot H^s_x(\R^2)} \lesssim A+$$
and hence (on undoing the scaling)
$$ \| u \|_{\dot H^s_x(\R^2)} \lesssim A \lambda^s$$
which (together with mass conservation) gives \eqref{apriori} as desired.
\end{proof}

\begin{remark}  By pursuing the above analysis more carefully, we in fact obtain a bound of the form
$$ \| u(T) \|_{H^s_x(\R^2)} \lesssim (1 + \|u_0\|_{H^s_x(\R^2)})^{C_s} (1+T)^{\frac{ s(1-s)}{2(2s-1)}+}$$ 
for some $C_s > 0$.  
\end{remark}

\begin{remark}
  \label{focusingremark}
Theorems \ref{eacq} and \ref{main} also hold for the focusing analog of
\eqref{eq:NLS} (replacing $|u|^2 u$ by $- |u|^2 u$) provided we also
assume $\|u_0 \|_{L^2_x(\R^2)} < \| Q \|_{L^2_x(\R^2)}$. Here $Q$ is the ground state
profile which arises as the unique(up to translations) positive
solution of $-Q + \Delta Q = - Q^3$. Indeed, most of the argument
remains unchanged (particularly those involving the local theory, or
multilinear estimates).  The only new difficulty arises when trying to use the energy $E(u)$ to control the kinetic component $\|u\|_{\dot H^1(\R^2)}^2$, since the potential energy component of the energy is now negative.  However, the sharp Gagliardo-Nirenberg inequality \cite{weinstein} allows one to achieve this control (losing a constant, of course) provided that $\|u\|_{L^2_x(\R^2)} < \|Q\|_{L^2_x(\R^2)}$, allowing one to continue the argument without difficulty. As the modifications are rather standard we do not detail them further here.
\end{remark}

It remains to prove Theorem \ref{eacq}.  There are clearly three components to this task: firstly,
to construct the functional $\tilde E$; secondly, to establish the fixed-time bound \eqref{fixed-bound}; and thirdly, to obtain the almost conservation law \eqref{acl}.  The first two tasks are straightforward and will be accomplished in the next two sections.  The third is substantially more difficult and will occupy the remainder of the paper.

\section{Construction of the modified energy functional}

We begin with the construction of the modified energy functional $\tilde E$.  As in previous literature on the $I$-method (e.g. \cite{CKSTT1}, \cite{CKSTT3}, \cite{CKSTT2}), it is convenient to introduce some notation for multilinear expressions involving $u$.

Let $k$ be an integer, let $\Sigma_k \subset (\R^2)^k$ denote the space
$$ \Sigma_k := \{ (\xi_1,\ldots,\xi_k) \in (\R^2)^k: \xi_1 + \ldots + \xi_k = 0 \},$$
with the measure induced from Lebesgue measure $d\xi_1 \ldots d\xi_{k-1}$ by pushing forward under the map 
$$(\xi_1,\ldots,\xi_{k-1}) \mapsto (\xi_1,\ldots,\xi_{k-1},-\xi_1-\ldots-\xi_{k-1}).$$
If $M: \Sigma_k \to \C$ is a smooth tempered symbol, and $u_1,\ldots,u_k \in \Schwartz(\R^2)$ are Schwartz functions, we define the $k$-linear functional
$$ \Lambda_k( M; u_1,\ldots,u_k ) :=
\frac{1}{(2\pi)^{2(k-1)}} \Re \int_{\Sigma_k} M(\xi_1,\ldots,\xi_k) \widehat{u_1}(\xi_1) \ldots \widehat{u_k}(\xi_k).$$
When $k$ is even, we abbreviate
$$ \Lambda_k( M; u ) := \Lambda_k( M; u, \overline{u}, \ldots, u, \overline{u} ).$$
We observe that the quantity $\Lambda_k(M;u)$ is invariant if one permutes the even arguments
$\xi_2, \xi_4,\ldots,\xi_k$ of $M$, the odd arguments $\xi_1,\xi_3,\ldots,\xi_{k-1}$ of $M$, as well as the additional symmetry
$$ M(\xi_1,\xi_2,\ldots,\xi_{k-1},\xi_k) \mapsto \overline{M}(\xi_2,\xi_1,\ldots,\xi_k,\xi_{k-1})$$
which swaps the odd and even arguments, and also conjugates $M$.  This generates a finite group $G_k$ of order $|G_k| = (k/2)! \times (k/2)! \times 2$ of symmetries, acting on $\Sigma_k$ and thus on the class $m$ of symbols.  This leads to the symmetrization rule
\begin{equation}\label{msym}
\Lambda_k( M; u ) = \Lambda_k( [M]_\sym; u)
\end{equation}
where $[M]_\sym := \frac{1}{|G_k|} \sum_{g \in G_k} gM$ is the $G_k$-symmetric component of $M$.

Using the above notation and the Fourier inversion formula, we observe that
$$ E(Iu) = \Lambda_2( \sigma_2 ; u ) + \Lambda_4( \sigma_4 ; u )$$
where
$$ \sigma_2(\xi_1,\xi_2) := - \frac{1}{2} \xi_1 m_1 \cdot \xi_2 m_2 = \frac{1}{2} |\xi_1|^2 m_1^2$$
and
$$\sigma_4(\xi_1,\xi_2,\xi_3,\xi_4) := \frac{1}{4} m_1 m_2 m_3 m_4 $$
and we abbreviate $m(\xi_j)$ as $m_j$.  Observe that $\sigma_2$ and $\sigma_4$ are both symmetric with respect to the group $G_k$.

Now we investigate the behaviour of these multilinear forms in time.  If $u$ is a smooth-in-time, Schwartz-in-space solution to \eqref{eq:NLS}, and $M$ is independent of time and symmetric with respect to $G_k$, then from the identity
$$ u_t = i \Delta u - i u \overline{u} u$$
arising from \eqref{eq:NLS}, together with some Fourier analysis, we have the differentiation formula
\begin{align*}
\partial_t \Lambda_k( M; u(t) ) &=
\Lambda_k( i M \alpha_k; u(t) ) - \Lambda_{k+2}( i k X(M); u(t) ) \\
&=
\Lambda_k( i M \alpha_k; u(t) ) - \Lambda_{k+2}( [i k X(M)]_\sym; u(t) )
\end{align*}
where $\alpha_k$ is the symbol
$$ \alpha_k(\xi_1,\ldots,\xi_k) := -|\xi_1|^2 + |\xi_2|^2 - \ldots - |\xi_{k-1}|^2 + |\xi_k|^2$$
(in particular, we have $\alpha_2 = 0$ on $\Sigma_2$) and $X(M)$ is the extended symbol
$$ X(M)(\xi_1,\ldots,\xi_{k+2}) := M( \xi_{123}, \xi_4, \ldots, \xi_{k+2})$$
where we use the notational convention $\xi_{ab} := \xi_a + \xi_b$, $\xi_{abc} := \xi_a + \xi_b + \xi_c$, etc.  Note that $iM\alpha_k$ is already symmetric with respect to $G_k$ and thus does not require further symmetrising.


As one particular instance of the above computations and the fundamental theorem of calculus, we have
\begin{align*}
E(Iu(t)) - E(Iu(0)) &= \int_0^t \partial_t E(Iu(t'))\ dt' \\
&= \int_0^t \Lambda_4( [- 2 i X(\sigma_2)]_\sym + i \sigma_4 \alpha_4; u(t') )\ dt'\\
&\quad - \int_0^t \Lambda_6( [4i X(\sigma_4)]_\sym; u(t') )\ dt'.
\end{align*}
In the case $m \equiv 1$ (which corresponds to $s=1$ or $N=\infty$), one easily computes that
$[-2iX(\sigma_2)]_\sym + i\sigma_4 \alpha_4$ and $[4iX(\sigma_4)]_\sym$ both vanish, thus giving a proof of energy conservation.  When $m$ is the multiplier from the previous section, these symbols do not vanish at high frequencies (when $\max(|\xi_1|,\ldots,|\xi_k|) \geq N$) but it turns out that the right-hand side can still be estimated by an expression which decays in $N$ as $O(N^{-3/2+})$; see \cite{CKSTT2}.  In fact, only the $\Lambda_4$ terms are as large as $O(N^{-3/2+})$; a closer inspection of the arguments in \cite{CKSTT2} show that the $\Lambda_6$ term is as least as small as $O(N^{-2+})$.  The strategy is thus to modify the quantity $E(Iu)$ so that the time derivative has less of a $\Lambda_4$ term and more of a $\Lambda_6$ term.  Specifically, we shall define
\begin{equation}\label{tedef}
\tilde E(u) := \Lambda_2( \sigma_2 ; u ) + \Lambda_4( \tilde \sigma_4 ; u )
\end{equation}
for some $G_4$-symmetric $\tilde \sigma_4$ to be chosen shortly.  Computing as before we have
\begin{align*}
\tilde E(u(t)) - \tilde E(u(0)) 
&= \int_0^t \Lambda_4( [- 2 i X(\sigma_2)]_\sym + i \tilde \sigma_4 \alpha_4; u(t') )\ dt'\\
&\quad - \int_0^t \Lambda_6( [4i X(\tilde \sigma_4)]_\sym; u(t') )\ dt'.
\end{align*}
An initial guess for $\tilde \sigma_4$ would thus be
$$ \tilde \sigma_4 := \frac{[2 i X(\sigma_2)]_\sym}{i \alpha_4}.$$
However this choice runs into the problem that $\alpha_4$ can vanish in the \emph{resonant interaction case} when $\xi_{12}$ and $\xi_{14}$ are either zero or orthogonal.  The first situation is easier to handle. In fact one can write 
\begin{equation}\label{alpha4}
\alpha_4 := -2 \xi_{12} \cdot \xi_{14} = -2 |\xi_{12}| |\xi_{14}| \cos \angle(\xi_{12},\xi_{14}) 
\end{equation}
and
\begin{equation}\label{sigma-20}
[2 i X(\sigma_2)]_\sym = \frac{1}{4} ( - m_1^2 |\xi_1|^2 + m_2^2 |\xi_2|^2 - m_3^2 |\xi_3|^2 + m_4^2 |\xi_4|^2 ). 
\end{equation}
In particular, when all frequencies are less than $N$ in magnitude, thus $\max_{1 \leq i \leq 4} |\xi_i| \leq N$, then we have computed
\begin{equation}\label{sigma-2}
[2 i X(\sigma_2)]_\sym = \frac{1}{4} i \alpha_4
\end{equation}
and so the vanishing of the denominator is cancelled by the numerator.
A similar argument can be used when $\xi_{12}=0$ or $
\xi_{14}=0$. Unfortunately, this cancellation is lost when one has one
or more high frequencies; this is in contrast to the one-dimensional
situation in \cite{CKSTT1}, \cite{CKSTT2}, where the resonant
interactions are simpler (and in \cite{CKSTT1}, one also has complete
integrability to provide further cancellations). 

Motivated by the above discussion, we shall in fact set
\begin{equation}\label{ts-def}
 \tilde \sigma_4 := \frac{[2 i X(\sigma_2)]_\sym}{i \alpha_4} 1_{\Omega_{nr}}
\end{equation}
where $1_{\Omega_{nr}}$ is the indicator function to the non-resonant set
\begin{equation}\label{omega-def}
\Omega_{nr} := \{ (\xi_1,\xi_2,\xi_3,\xi_4) \in \Sigma_4:
\max_{1 \leq j \leq 4} |\xi_j| \leq N \}
\cup \{ (\xi_1,\xi_2,\xi_3,\xi_4) \in \Sigma_4: |\cos \angle(\xi_{12}, \xi_{14})| \geq \theta_0 \},
\end{equation}
where $0 < \theta_0 < 1/100$ is a parameter to be chosen later (we will shortly take $\theta_0 := 1/N$).  Note that while the angle $\angle(\xi_{12}, \xi_{14})$ is undefined when $\xi_{12}$ or $\xi_{14}$ vanishes, but this set has measure zero and can be ignored.

\begin{remark}  The presence of the expression $|\cos \angle(\xi_{12}, \xi_{14})|$ is the key to all of our improvements over the previous work in \cite{CKSTT2}.  However, as this expression involves three of the four frequencies in $\Sigma_4$, exploiting this term properly will turn out to be a significant technical headache, requiring many decompositions of the frequency variables to handle.
\end{remark}

We now define $\tilde E$ by \eqref{tedef} with $\tilde{\sigma}_4$ as
in \eqref{ts-def}.  To prove Theorem \ref{eacq}, it thus suffices to prove the following two propositions.

\begin{proposition}[Fixed-time estimate]\label{fixed}  Let the notation be as above.  Then for any $u \in \Schwartz_x(\R^2)$, we have
\begin{equation}\label{fixed-bound2}
 | \Lambda_4( \sigma_4 - \tilde \sigma_4; u ) | \lesssim \theta_0^{-1} N^{-2+} \| Iu \|_{H^1_x(\R^2)}^4.
\end{equation}
\end{proposition}

\begin{proposition}[Almost conservation law]\label{acprop}  Let the notation be as above.  
If $\|u_0\|_{L^2_x(\R^2)} \leq A$ and $E(Iu_0) \leq 1$, and $u$ is a smooth-in-time, Schwartz-in-space solution to \eqref{eq:NLS} on a time interval $[0,t_0]$, then if $t_0$ is sufficiently small depending on $A$, we have
$$
|\int_0^{t_0} \Lambda_4( [- 2 i X(\sigma_2)]_\sym + i \tilde \sigma_4 \alpha_4; u(t) )\ dt|
\lesssim C(A) [N^{-2+} + \theta_0^{1/2} N^{-3/2+} + \theta_0^{-1} N^{-3+}]$$
and
$$ |\int_0^{t_0} \Lambda_6( [4i X(\tilde \sigma_4)]_\sym; u(t) )\ dt| \lesssim C(A) [N^{-2+} + \theta_0^{1/2} N^{-3/2+} + \theta_0^{-1} N^{-3+}].$$
\end{proposition}

Indeed, by setting $\theta_0 := 1/N$ we obtain the desired result. In
fact, we will see below (see Remark \ref{teetertotter} below and the
two propositions preceding it) that the 6-linear estimate degenerates
with growing $\theta_0$ while the 4-linear estimate improves with
$\theta_0$ and that the choice $\theta_0 = 1/N$ puts these contributions to
the energy increment in balance.

The rest of the paper is now devoted to the proof of these two propositions.

\section{The fixed time estimate}

In this section we prove Proposition \ref{fixed}, which is in fact rather easy.  
From Plancherel's theorem, it suffices to show that
$$ \int_{\Sigma_4} |\sigma_4(\xi) - \tilde \sigma_4(\xi)| \prod_{j=1}^4 \frac{|\hat u_j(\xi_j)|}{m_j \langle \xi_j \rangle}
\lesssim \theta_0^{-1} N^{-2+} \prod_{j=1}^4 \|u_j\|_{L^2(\R^2)}$$
for any $u_1,u_2,u_3,u_4 \in L^2(\R^2)$, where $\xi := (\xi_1,\xi_2,\xi_3,\xi_4)$.

From \eqref{sigma-2}, \eqref{ts-def}, \eqref{omega-def} we know that $\sigma_4(\xi)-\tilde \sigma_4(\xi)$ vanishes when $\max_{1\leq j\leq 4} |\xi_j| \leq N$, so we may restrict to
the region $\max_{1\leq j\leq 4} |\xi_j| > N$.  We now need the following bound.

\begin{lemma}\label{symbo}  For any $(\xi_1,\xi_2,\xi_3,\xi_4) \in \Sigma_4$, We have
$$ |[2 i X(\sigma_2)]_\sym| \lesssim \min( m_1,m_2,m_3,m_4)^2 |\xi_{12}| |\xi_{14}|.$$
\end{lemma}

\begin{proof} Let $f(\xi) := m(\xi)^2 |\xi|^2$.  In light of \eqref{sigma-20}, it suffices
to show that
$$
 |f(\xi_1) - f(\xi_2) + f(\xi_3) - f(\xi_4)| \lesssim \min( m_1,m_2,m_3,m_4)^2 |\xi_{12}| |\xi_{14}|.
$$
Using symmetries, we may assume that $|\xi_1| \geq |\xi_2|, |\xi_3|, |\xi_4|$ and
$|\xi_{12}| \geq |\xi_{14}|$.  In particular $\min(m_1,m_2,m_3,m_4) = m_1$.

First assume that $|\xi_{12}|, |\xi_{14}| \gtrsim |\xi_1|$.  Then we can estimate all four terms on the left-hand side by $O( m_1^2 |\xi_1|^2 )$, and the claim follows.

Now assume that $|\xi_{12}| \sim |\xi_1|$ but that $|\xi_{14}| \ll |\xi_1|$.  We write the left-hand side as
$$ | (f(\xi_1) - f(\xi_1 - f(\xi_{14}))) + (f(\xi_3) - f(\xi_3 + \xi_{14})) |.$$
Note that $\nabla f(\xi) = O( m(\xi)^2 |\xi| )$, and $m(\xi)^2 |\xi|$ is an increasing function of $|\xi|$, so by the fundamental theorem of calculus we have
$$ |f(\xi_1) - f(\xi_1 - f(\xi_{14}))|, |f(\xi_3) - f(\xi_3 + \xi_{14})| \lesssim m_1^2 |\xi_1| |\xi_{14}|$$
and the claim follows.

Finally, suppose that $|\xi_{12}|, |\xi_{14}| \ll |\xi_1|$.
We write the left-hand side as
$$ |f(\xi_1) - f(\xi_1 - \xi_{12}) - f(\xi_1 - \xi_{14}) + f(\xi_1 - \xi_{12} - \xi_{14})|$$
which we can write as
$$ |\int_0^1 \int_0^1 (\xi_{12} \cdot \nabla) (\xi_{14} \cdot \nabla) f( \xi_1 - s \xi_{12} - t \xi_{14})\ ds dt|.$$
Since $\nabla^2 f(\xi_1-s\xi_{12}-t\xi_{14})| = O(m_1^2)$, the claim follows.
\end{proof}

From this lemma and \eqref{ts-def}, \eqref{omega-def}, \eqref{alpha4}, we obtain the following useful pointwise bound:

\begin{corollary}\label{ts-bound}  For any $(\xi_1, \xi_2, \xi_3,
  \xi_4) \in \Sigma_4$, we have
$$ |\tilde \sigma_4| \lesssim \frac{\min(m_1,m_2,m_3,m_4)^2}{\theta_0}.$$
\end{corollary}

Since, for $(\xi_1, \xi_2, \xi_3,
  \xi_4) \in \Sigma_4$, we have
$$|\sigma_4| \sim m_1 m_2 m_3 m_4 \lesssim \min(m_1,m_2,m_3,m_4)^2,$$
we reduce matters to showing that
$$ \int_{\Sigma_4: \max_{1\leq j\leq 4} |\xi_j| > N} \min(m_1,m_2,m_3,m_4)^2 \prod_{j=1}^4 \frac{|\hat u_j(\xi_j)|}{m_j \langle \xi_j \rangle}
\lesssim N^{-2} \prod_{j=1}^4 \|u_j\|_{L^2(\R^2)}.$$
Note that at least two of $|\xi_1|, |\xi_2|, |\xi_3|, |\xi_4|$ need to be greater than or comparable to $N$.
Without loss of generality we may assume that $|\xi_1|, |\xi_2| \gtrsim N$.  Then
$m_j \langle \xi_j \rangle \gtrsim N^{1-} \langle \xi_j \rangle^{0+}$ for $j=1,2$, while
$\min(m_1,m_2,m_3,m_4)^2 \lesssim m_3 m_4$, so we reduce to showing that
$$ \int_{\Sigma_4} \langle \xi_1 \rangle^{0-} \langle \xi_2 \rangle^{0-}
\langle \xi_3 \rangle^{-1} \langle \xi_4 \rangle^{-1} \prod_{j=1}^4 |\hat u_j(\xi_j)|
\lesssim \prod_{j=1}^4 \|u_j\|_{L^2(\R^2)},$$
which by Plancherel is equivalent to the estimate
$$ |\int_{\R^2} v_1 v_2 v_3 v_4\ dx| \lesssim \| v_1 \|_{H^{0+}(\R^2)}
\| v_2 \|_{H^{0+}(\R^2)} \| v_3 \|_{H^1(\R^2)} \| v_4 \|_{H^1(\R^2)}
$$
for some $v_1,v_2,v_3,v_4$.  But this easily follows from Sobolev embedding and H\"older.  This proves
Proposition \ref{fixed}. Note that this proof only required $s >0$.

\section{Modified local well-posedness}

It remains to prove Proposition \ref{acprop}.  From the hypotheses on $u_0$ we have
$$ \| Iu_0 \|_{H^1_x(\R^2)} \lesssim A.$$
In order to use this bound, we need some spacetime estimates on the solution $u$.  We recall
the standard $X^{s,b}(\R \times \R^2)$ spaces for $s,b \in \R$, defined on spacetime Schwartz functions by the norm\footnote{Note that our sign conventions for the Schr\"odinger equation and the spacetime Fourier transform force the dispersion relation to be $\tau = -|\xi|^2$ rather than $\tau = +|\xi|^2$.  Of course, these sign conventions are not crucial to our final results.}
$$ \| u \|_{X^{s,b}(\R \times \R^2)} := \| \langle \xi \rangle^s 
\langle \tau + |\xi|^2 \rangle^b \tilde u(\tau,\xi) \|_{L^2_{\tau,\xi}(\R \times \R^2)},$$
where
$$ \tilde u(\tau,\xi) := \int_\R \int_{\R^2} e^{-i (t\tau +x \cdot \xi)} u(t,x)\ dx dt$$
is the spacetime Fourier transform of $u$, and then for any time interval $J$, define the restricted
norm $X^{s,b}(J \times \R^2)$ by the formula
$$\|u\|_{X^{s,b}(J \times \R^2)} := \inf \{ \|v\|_{X^{s,b}(\R \times \R^2)}: v|_{J \times \R^2} = u \}$$
where $v$ ranges over all functions in $X^{s,b}(\R \times \R^2)$ which agree with $u$ on $J \times \R^2$.  We caution that $u$ and $\overline{u}$ need not have comparable $X^{s,b}$ norms; this will complicate our notation a little bit but will not significantly affect the analysis.

We now fix an exponent $b$ close to $1/2$ (e.g. $b := 0.6$).

\begin{proposition}[Modified local existence]\label{mle}  Let $u_0$ be such that $\|Iu_0\|_{H^1_x(\R^2)} \lesssim A$, and $u$ is a smooth-in-time, Schwartz-in-space solution to \eqref{eq:NLS} on a time interval $[0,t_0]$, then if $t_0$ is sufficiently small depending on $A$, we have
$$ \| Iu \|_{X^{1,b}([0,t_0] \times \R^2)} \lesssim A.$$
\end{proposition}

\begin{proof}  See \cite[Proposition 3.2]{CKSTT2}.  The proposition
  there was stated only for $s > 4/7$ and for an unspecified $b$, but
  it is not difficult to see that the argument in fact works for $b :=
  0.6$ and for all $s > 1/2$.  (In fact, the argument works for all $s > 0$, though as $s$ approaches $0$ one needs to let $b$ approach $1/2$.) 
\end{proof}

In view of this proposition, we see that to prove Proposition \ref{acprop} it suffices to prove the following estimates.

\begin{proposition}[Quadrilinear estimate]\label{quatro}  For any $Iu \in X^{1,b}(\R \times \R^2)$ and $0 < t_0 < 1$ we have
\begin{equation}\label{4-est}
|\int_0^{t_0} \Lambda_4( [- 2 i X(\sigma_2)]_\sym + i \tilde \sigma_4 \alpha_4; u(t) )\ dt|
\lesssim [N^{-2+} + \theta_0^{1/2} N^{-3/2+}] \| Iu \|_{X^{1,b}(\R \times \R^2)}^4.
\end{equation}
\end{proposition}

\begin{proposition}[Sextilinear estimate]\label{sex} For any $Iu \in X^{1,b}(\R \times \R^2)$ and $0 < t_0 < 1$ we have
\begin{equation}\label{6-est}
 |\int_0^{t_0} \Lambda_6( [4i X(\tilde \sigma_4)]_\sym; u(t) )\ dt| \lesssim \theta_0^{-1} N^{-3+} \| Iu \|_{X^{1,b}(\R \times \R^2)}^6.
\end{equation}
\end{proposition}

\begin{remark} 
\label{teetertotter} Observe that decreasing the threshold $\theta_0$ between resonance and non-resonance improves the quadrilinear estimate (fewer resonant interactions) at the expense of the sextilinear estimate (more non-resonant interactions).  The case $\theta=1$ is essentially the case considered in \cite{CKSTT2}.
\end{remark}

The proof of these propositions will occupy the remainder of the paper.

Henceforth all spacetime norms will be on the full spacetime domain $\R \times \R^2$, and we shall omit this domain from the notation for brevity.

\section{$X^{s,b}$ estimates}

In this section we record some standard estimates involving the $X^{s,b}$ spaces which we will need in the sequel.

Let us say that a function $u$ has \emph{spatial frequency $N$} if its Fourier transform (either spatial or spacetime) is supported on the annulus $\{ \langle \xi\rangle \sim N \}$.  From the standard energy estimate $\|u\|_{L^\infty_t L^2_x} \lesssim \|u\|_{X^{0,1/2+}}$ (see e.g. \cite[Corollary 2.10]{tao}) and Bernstein's inequality we have

\begin{lemma}[Linear estimate]  
If $u$ has spatial frequency $N$, then
\begin{equation}\label{sup-est}
\|u\|_{L^\infty_{t,x}} \lesssim N^{1} \| u \|_{X^{0,1/2+}}.
\end{equation}
\end{lemma}

We also have some standard bilinear estimates:

\begin{lemma}[Bilinear Strichartz estimate]  
If $u_1$, $u_2$ have spatial frequency $N_1$, $N_2$ respectively, then
\begin{equation}\label{bil-est}
\|u_1 u_2 \|_{L^2_{t,x}} \lesssim \frac{N_2^{1/2}}{N_1^{1/2}}
 \| u_1 \|_{X^{0,1/2+}}  \| u_2 \|_{X^{0,1/2+}}.
\end{equation}
If, furthermore, $N_2 \ll N_1$ and $u_2$ has Fourier support supported in a ball\footnote{This ball is actually a cylinder if one also considers the time-frequency variable $\tau$.} $\{ \xi = \xi_0 + O(\theta N_2) \}$ of radius $O(\theta N_2)$ for some $0 < \theta < 1$, then we can improve the above estimate to
\begin{equation}\label{bil-est-theta}
\|u_1 u_2 \|_{L^2_{t,x}} \lesssim \frac{\theta^{1/2} N_2^{1/2}}{N_1^{1/2}}
 \| u_1 \|_{X^{0,1/2+}}  \| u_2 \|_{X^{0,1/2+}}.
\end{equation}
\end{lemma}

\begin{proof}  The estimate \eqref{bil-est} is standard, see e.g. \cite{Bo1}, \cite{Bo2} (see also \cite[Lemma 3.4]{gopher}).  The second claim then follows by a Galilean transformation argument, shifting the frequencies of $u_1, u_2$ by about $N_2$ to ensure that $u_2$ now has frequency $\sim \theta N_2$ rather than $N_2$, without significantly affecting the frequency of $N_1$.  Note that Galilean transforms do not affect the $X^{0,1/2+}$ norm or the $L^2_{t,x}$ norm of $u_1 u_2$.
\end{proof}

Of course, it is advantageous to apply this estimate when $N_2 \leq N_1$ rather than when $N_1 \leq N_2$.  We also make the trivial remark that we can replace $u_1 u_2$ by $\overline{u_1} u_2$, $u_1 \overline{u_2}$, or $\overline{u_1} \overline{u_2}$ without affecting the estimate.

\section{Proof of sextilinear estimate}


In this section we prove Proposition \ref{sex}, which is the easier of the two propositions, as it does not require any fine control on the resonant interactions\footnote{Indeed, there seems to be a general principle when applying the $I$-method that terms which are more multilinear (and hence have fewer derivatives) are easier to estimate than terms which are less multilinear.  This phenomenon, which is related to the sub-criticality of the regularities being considered, explains why it is beneficial to adjust the modified energy using correction terms, as this makes the error terms more multilinear.}. 
The left-hand side of \eqref{6-est} can be expanded as
$$ |\int_0^{t_0} [ \int_{\Sigma_6} [4i X(\tilde \sigma_4)]_\sym(\xi)
\hat u(t,\xi_1) \ldots \hat{\overline{u}}(t,\xi_6) ]\ dt|.$$
If $\max(|\xi_1|,\ldots,|\xi_6|) < N/3$, then (by \eqref{sigma-2}, \eqref{ts-def}), $4iX(\tilde \sigma_4) = 1$, and thus $[4i X(\tilde \sigma_4)]_\sym$
vanishes.  We can thus restrict to the region $\max(\xi_1,\ldots,\xi_6) \geq N/3$.  We then remove the symmetry and reduce to showing that
$$ |\int_0^{t_0} [ \int_{\Sigma_6: \max(|\xi_1|,\ldots,|\xi_6|) \geq N/3} X(\tilde \sigma_4)
\hat u(t,\xi_1) \ldots \hat{\overline{u}}(t,\xi_6) ]\ dt| \lesssim \theta_0^{-1} N^{-3+}
\|Iu\|_{X^{1,b}}^6.$$
Because the $X^{s,b}$ norm uses the spacetime Fourier transform, we will be forced for technical reasons\footnote{The specific issue is that we cannot automatically reduce to the case where the spatial Fourier transforms of the $u_j$ are non-negative.  In previous literature this difficulty was avoided by using the Coifman-Meyer multiplier theorem, but the symbol here does not obey Coifman-Meyer type estimates and so this theorem is not applicable.} to write the left-hand side 
in terms of the spacetime Fourier transform.  Indeed, this left-hand side becomes
$$ |\int_\R \ldots \int_\R [\int_{\Sigma_6} 
\hat 1_{[0,t_0]}(\tau_0) 1_{\max(|\xi_1|,\ldots,|\xi_6|) \geq N/3} X(\tilde \sigma_4)\tilde u(\tau_1,\xi_1)
\ldots \tilde {\overline{u}}(\tau_6,\xi_6)] d\tau_1 \ldots d\tau_6|$$
where $\tau_0 := -\tau_1 - \ldots -\tau_6$.  Using Corollary \ref{ts-bound} and the bound
\begin{equation}\label{t0-fourier}
\hat 1_{[0,t_0]}(\tau_0) = \langle \tau_0 \rangle^{-1},
\end{equation}
we can estimate this quantity by
$$ \lesssim \frac{1}{\theta_0}
\int\limits_\R \ldots \int\limits_\R [\int\limits_{\Sigma_6} 
\langle \tau_0 \rangle^{-1} 1_{\max(|\xi_1|,\ldots,|\xi_6|) \geq N/3} 
\min( m_{123}, m_4, m_5, m_6 )^2 |\tilde u(\tau_1,\xi_1)|
\ldots |\tilde{\overline{u}}(\tau_6,\xi_6)|] d\tau_1 \ldots d\tau_6$$
where $m_{123} = m(\xi_{123})$.  It will be convenient to hide the conjugations by using the norm
$$ \| u \|_{\tilde X^{s,b}} := \inf \{ \| u_1 \|_{X^{s,b}} + \| \overline{u_2} \|_{X^{s,b}}: u = u_1 + u_2 \}$$
so it suffices to show the estimate
\begin{eqnarray*} \int\limits_\R  \ldots \int\limits_\R & [\int\limits_{\Sigma_6} 
\langle \tau_0 \rangle^{-1} 1_{\max(|\xi_1|,\ldots,|\xi_6|) \geq N/3} 
\min( m_{123}, m_4, m_5, m_6 )^2 \prod_{j=1}^6 |\tilde u_j(\tau_j,\xi_j)|] 
d\tau_1 \ldots d\tau_6\\
& \lesssim N^{-3+} \prod_{j=1}^6 \|Iu_j\|_{\tilde X^{1,b}}.\end{eqnarray*}
If we let $|\xi^*_1| \geq \ldots \geq |\xi^*_6|$ be the six magnitudes $|\xi_1|,\ldots,|\xi_6|$
in order, we observe that $\min( m_{123}, m_4, m_5, m_6 )^2 \leq
m(\xi^*_4)^2$.  Inserting this bound, the left-hand side is now
symmetric 
in $\xi_1,\ldots,\xi_6$, so we can assume that
$|\xi_1| \geq \ldots \geq |\xi_6|$.  The constraint $\max(|\xi_1|,\ldots,|\xi_6|) \geq N/3$ then implies $|\xi_2| \gtrsim N$.  We thus need to show that
$$ \int_\R \ldots \int_\R [\int_{\Sigma_6} 
\langle \tau_0 \rangle^{-1} 1_{|\xi_1| \geq \ldots \geq |\xi_6|} 1_{|\xi_2| \gtrsim N}
m(\xi_4)^2 \prod_{j=1}^6 |\tilde u_j(\tau_j,\xi_j)|] 
d\tau_1 \ldots d\tau_6 \lesssim N^{-3+} \prod_{j=1}^6 \|Iu_j\|_{\tilde X^{1,b}}.$$
Partitioning up into Littlewood-Paley pieces, it suffices to show that
$$ \int_\R \ldots \int_\R [\int_{\Sigma_6} 
\langle \tau_0 \rangle^{-1} \prod_{j=1}^6 |\tilde u_j(\tau_j,\xi_j)|] 
d\tau_1 \ldots d\tau_6 \lesssim \frac{N^{-3+} N_1^{0-} \prod_{j=1}^6 m(N_j) N_j \|u_j\|_{\tilde X^{0,b}}}{m(N_4)^2} $$
whenever $N_1 \gtrsim N_2 \gtrsim \ldots \gtrsim N_6 \gtrsim 1$, and each $u_j$ has spatial frequency $N_j$.  Note we may assume $N_1 \sim N_2$ since the $\Sigma_6$ integral vanishes otherwise.
Since the definition of $\tilde X^{1,b}$ uses only the magnitude of the spacetime Fourier transform, we may take all of the $\tilde u_j$ to be non-negative and thus omit the absolute value signs.
The left-hand side can now be written using spacetime convolutions as
\begin{equation}\label{mn}
 \int_\R \langle \tau \rangle^{-1} \tilde u_1 * \ldots * \tilde u_6(\tau,0)\ d\tau.
\end{equation}
It is slightly unfortunate that $\langle \tau \rangle^{-1}$ barely fails to be integrable.  However, if we introduce the logarithmic weight $w(\tau) := 1 + \log^2 \langle \tau \rangle$, then $\langle \tau \rangle^{-1} w^{-1}$ is integrable.  Also, from the elementary estimate $w(\tau_1 + \ldots + \tau_6) \lesssim w(\tau_1) \ldots w(\tau_6)$ we have the pointwise bound
$$ \tilde u_1 * \ldots * \tilde u_6 \lesssim w^{-1} [ (w \tilde u_1) * \ldots (w \tilde u_6) ]$$
and thus we can bound
$$ \eqref{mn} \lesssim \| (w\tilde u_1) * \ldots * (w\tilde u_6) \|_{L^\infty_{\tau,\xi}}.$$
Thus it will suffice to show that
$$\| (w\tilde u_1) * \ldots * (w\tilde u_6) \|_{L^\infty_{\tau,\xi}}
\lesssim \frac{N^{-3+} N_1^{0-} \prod_{j=1}^6 m(N_j) N_j \| u_j \|_{\tilde X^{0,b}}}{m(N_4)^2}.$$
If $v_j$ denotes the function with spacetime Fourier transform $\tilde v_j = w \tilde u_j$, one easily verifies that
$$ \| v_j \|_{\tilde X^{0,b-}} \lesssim \log^2 (1 + N_1) \| u_j \|_{\tilde X^{0,b}}$$
and so it will suffice to show that
$$ \| \tilde v_1 * \ldots * \tilde v_6 \|_{L^\infty_{\tau,\xi}}
\lesssim \frac{N^{-3+} N_1^{0-} \prod_{j=1}^6 m(N_j) N_j \| v_j \|_{\tilde X^{0,b-}}}{m(N_4)^2}.$$
By Hausdorff-Young it suffices to show that
$$ \| v_1 \ldots v_6 \|_{L^1_{t,x}(\R \times \R^2)}
\lesssim \frac{N^{-3+} N_1^{0-} \prod_{j=1}^6 m(N_j) N_j \| v_j \|_{\tilde X^{0,b-}}}{m(N_4)^2}.$$
Since the left-hand side is insensitive to conjugation, it suffices to show that
$$ \| v_1 \ldots v_6 \|_{L^1_{t,x}(\R \times \R^2)}
\lesssim \frac{N^{-3+} N_1^{0-} \prod_{j=1}^6 m(N_j) N_j \| v_j \|_{X^{0,b-}(\R \times \R^2)}}{m(N_4)^2}.$$
Estimating $v_1 v_3$ and $v_2 v_4$ in $L^2_{t,x}$ using
\eqref{bil-est}, and $v_5$, $v_6$ in $L^\infty_{t,x}$ using
\eqref{sup-est}, and applying H\"older's inequality, we reduce to
showing that
%
%
%
$$ \frac{N_3^{1/2}}{N_1^{1/2}} \frac{N_4^{1/2}}{N_2^{1/2}} N_5^{1} N_6^{1}
 \lesssim \frac{N^{-3+} N_1^{0-} \prod_{j=1}^6 m(N_j) N_j}{m(N_4)^2}.$$
 We rearrange (using $N_1 \sim N_2$) as
 $$ 1 \lesssim N^{-3+} m(N_1)^2 N_1^{3-} m(N_3) N_3^{1/2} m(N_4) N_4^{1/2}
 m(N_4)^{-2} m(N_5)  m(N_6) .$$
 Since the function $m(|\xi|)$ is non-increasing in $|\xi|$, we may reduce to the case $N_4=N_5=N_6$,
 which becomes
 $$ 1 \lesssim N^{-3+} m(N_1)^2 N_1^{3-} m(N_3) N_3^{1/2} m(N_4) N_4^{1/2}.$$
This is true since $m(N_3) N_3^{1/2} \geq m(N_4) N_4^{1/2} \geq 1$ and $m(N_1) N_1
\gtrsim N$. This concludes the proof of Proposition \ref{sex}.

\section{Proof of quadrilinear estimate}

We now begin the proof of Proposition \ref{quatro}.  We shall begin by mimicking the proof of
Proposition \ref{sex}, but we will find that there are a few cases, particularly those involving resonant interactions, which require more careful attention, both in the pointwise estimates on the multiplier, and on the bilinear estimates needed to handle the final expression.

From \eqref{ts-def} we have
$$ ([- 2 i X(\sigma_2)]_\sym + i \tilde \sigma_4 \alpha_4)(\xi) = [- 2 i X(\sigma_2)]_\sym 1_{\Omega_r}$$
where the \emph{resonant set} 
$$ \Omega_r := \{ (\xi_1,\xi_2,\xi_3,\xi_4) \in \Sigma_4: \max(|\xi_1|,|\xi_2|,|\xi_3|,|\xi_4|) > N;
|\cos \angle(\xi_{12}, \xi_{14})| < \theta_0 \}$$
is the complement of $\Omega_{nr}$.  Thus the left-hand side of \eqref{4-est} can be expressed using
spacetime Fourier transforms similarly to the previous section as
$$ |\int_\R \ldots \int_\R \widehat{1_{[0,t_0]}}(\tau_0) [ \int_{\Omega_{r}} [- 2 i X(\sigma_2)]_\sym(\xi)
\tilde u(\tau_1,\xi_1) \ldots \tilde{\overline{u}}(\tau_4,\xi_4) ]\ d\tau_1 \ldots d\tau_4|$$
where $\tau_0 := -\tau_1 - \ldots - \tau_4$.
Using the bound \eqref{t0-fourier}, we can bound this by
$$ \lesssim
\int_\R \ldots \int_\R \langle \tau_0 \rangle^{-1} [\int_{\Omega_{r}}
|[X(\sigma_2)]_\sym|
|\tilde u(\tau_1,\xi_1)| \ldots |\tilde{\overline{u}}(\tau_4,\xi_4)| ]\ d\tau_1 \ldots d\tau_4.$$
By symmetry we may reduce to the region where $|\xi_1| \geq |\xi_2|, |\xi_3|, |\xi_4|$ and
$|\xi_{2}| \geq |\xi_{4}|$, or equivalently we may replace $\Omega_{r}$ with
$$ \Omega'_r := \{ (\xi_1,\xi_2,\xi_3,\xi_4) \in \Sigma_4: 
|\xi_1| \geq |\xi_2|, |\xi_3|, |\xi_4|; |\xi_{2}| \geq |\xi_{4}|; |\xi_1| > N;
|\cos \angle(\xi_{12}, \xi_{14})| < \theta_0 \}.$$
It thus suffices as in the previous section to show the estimate
\begin{align*}
&\int_\R \ldots \int_\R \langle \tau_0 \rangle^{-1} [\int_{\Omega'_r}
|[X(\sigma_2)]_\sym|| \prod_{j=1}^4 \tilde u^*_j(\tau_j,\xi_j) ]\ d\tau_1 \ldots d\tau_4\\
&\quad \lesssim N_1^{0-} [N^{-2+} + \theta_0^{1/2} N^{-3/2+}]  \prod_{j=1}^4
m(N_j) N_j \| u_j \|_{\tilde{X}^{0,b}(\R \times \R^2)}
\end{align*}
whenever 
$$
N_1 \gtrsim N, N_2,N_3,N_4 \gtrsim 1; \quad N_2 \gtrsim N_4,
$$
and $u_1,\ldots,u_4$ have spatial frequency $N_1,\ldots,N_4$ respectively, and have non-negative spacetime Fourier transform.  Here we adopt the convention that $u^*_j = u_j$ when $j$ is odd and $u^*_j = \overline{u_j}$ when $j$ is even.  Note that we may also assume that
$$
N_1 \sim N_2
$$
for if $N_2 \ll N_1$ then $\xi_{12}, \xi_{14}$ make a small angle with $\xi_1$ and it will be impossible to satisfy the condition $|\cos \angle(\xi_{12}, \xi_{14})| \leq \theta_0$.  The constraints on $N_1,N_2,N_3,N_4$ are now symmetric under swapping $N_1, N_3$ with $N_2,N_4$ respectively and so we may now also assume $N_3 \gtrsim N_4$.
To summarize, the frequencies $N_1, N_2, N_3, N_4$ are now known to obey the relations
\begin{equation}\label{n-relations}
N_1 \sim N_2 \gtrsim N; \quad N_2 \gtrsim N_3 \gtrsim N_4 \gtrsim 1.
\end{equation}
Using the weight $w$ and the functions $v_j$ as in the previous section, we reduce to showing that
\begin{equation}\label{nmonster}
\begin{split}
&\sup_{\tau} \int_\R \ldots \int_\R  [\int_{\Omega'_r}
|[X(\sigma_2)]_\sym| \prod_{j=1}^4 \tilde v^*_j(\tau_j,\xi_j) ]
\delta(\tau_0-\tau) d \tau_1 \dots d\tau_4 \\
&\quad \lesssim N_1^{0-} [N^{-2+} + \theta_0^{1/2} N^{-3/2+}]  \prod_{j=1}^4
m(N_j) N_j \| v_j \|_{\tilde{X}^{0,b-}}.
\end{split}
\end{equation}

We now dispose of an easy case, in which we will not use the $\theta_0^{1/2} N^{-3/2+}$ term on the right-hand side.  Suppose that
$$ N_4 \gtrsim N \hbox{ or } N_3 \sim N_4.$$
Using Lemma \ref{symbo} we have
\begin{equation}\label{xs}
[X(\sigma_2)]_\sym \lesssim m(N_1)^2 |\xi_{14}| |\xi_{12}| \lesssim m(N_1)^2 N_1 N_3
\end{equation}
since $|\xi_{12}| = |\xi_{34}|$ is bounded by $N_3$.  Gathering some
terms and simplifying using \eqref{n-relations}, we reduce to showing that
\begin{align*}
&\sup_{\tau} \int_\R \ldots \int_\R  [\int_{\Omega'_r}
\prod_{j=1}^4 \tilde v^*_j(\tau_j,\xi_j) ] \delta(\tau_0-\tau) d\tau_1
\cdots d\tau_4 \\
&\quad \lesssim N_1^{1-} [N^{-2+} + \theta_0^{1/2} N^{-3/2+}] m(N_3) m(N_4) N_4 \prod_{j=1}^4
\| v_j \|_{\tilde{X}^{0,b-}}.
\end{align*}
We estimate the left-hand side by
$$ \| \tilde v^*_1 * \ldots * \tilde v^*_4 \|_{L^\infty_{\tau,\xi}}
\lesssim \| v_1 \overline{v_2} v_3 \overline{v_4} \|_{L^1_{t,x}}.$$
Applying \eqref{bil-est} to $v_1 v_3$ and $v_2 v_4$, followed by Cauchy-Schwarz, we reduce to showing that
$$ \frac{N_3^{1/2}}{N_1^{1/2}} \frac{N_4^{1/2}}{N_2^{1/2}}
\lesssim N_1^{1-} N^{-2+} m(N_3) m(N_4) N_4$$
which simplifies using \eqref{n-relations} to
$$ 1 \lesssim N^{-2+} N_1^{2-} m(N_3) N_3^{-1/2} m(N_4) N_4^{1/2}.$$
If $N_3 \sim N_4$, then this becomes
$$ 1 \lesssim N^{-2+} N_1^{2-} m(N_3)^2;$$
since $m(N_1)^2 \lesssim m(N_3)^2$, $m(\xi)^{2-} |\xi|^2$ is essentially increasing in $|\xi|$, and $N_1 \gtrsim N$, we conclude the claim.  If instead $N_4 \gtrsim N$, we have $m(N_4) N_4^{1/2} \gtrsim N^{1/2}$ and $m(N_3) N_3^{-1/2} \gtrsim m(N_1) N_1^{-1/2}$, thus reducing to
$$ 1 \lesssim N^{-3/2+} N_1^{3/2-} m(N_1)$$
which follows since $N_1 \gtrsim N$.

To conclude Theorem \ref{main} we only need to establish
\eqref{nmonster} in the remaining case when
\begin{equation}\label{n-relations-2}
N_1 \sim N_2 \gtrsim N, N_3 \gg N_4 \gtrsim 1.
\end{equation}
It is here that we shall need to exploit the resonance constraint $|\cos \angle(\xi_{12},\xi_{14})| \leq \theta_0$ more fully.  Specifically, we shall need the following improved bound on $|[X(\sigma_2)]_\sym|$.

\begin{lemma}  Let $N_1,\ldots,N_4$ be as in \eqref{n-relations-2}, and let $(\xi_1,\xi_2,\xi_3,\xi_4) \in \Omega'_r$ be such that $|\xi_j| \sim N_j$ for $j=1,2,3,4$. Then
$$ |[X(\sigma_2)]_\sym| \lesssim m(N_1)^2 N_1 N_3 \theta_0 +  m(N_3)^2 N_3^2$$
(compare with \eqref{xs}).
\end{lemma}

\begin{proof} The new idea is to exploit heavily the spherical symmetry of $m$.
From \eqref{sigma-20} we have
$$ |[X(\sigma_2)]_\sym| \lesssim \left| m(|\xi_1|)^2 |\xi_1|^2 - m(|\xi_2|)^2 |\xi_2|^2 \right| + m(N_3)^2 N_3^2 + m(N_4)^2 N_4^2.$$
Since the function $x \mapsto m(x)^2 x^2$ is essentially increasing, it suffices to show that
$$ m(|\xi_1|)^2 |\xi_1|^2 - m(|\xi_2|)^2 |\xi_2|^2 = m(N_1)^2(O( N_3 N_1 \theta_0 ) + O(N_3)^2).$$
On the region $x \sim N_1$, the function $x \mapsto m(x)^2 x^2$ has derivative $O( m(N_1)^2 N_1 )$.  Thus we can
bound the left-hand side by
$$ m(N_1)^2 N_1 | |\xi_1| - |\xi_2| | \sim m(N_1)^2 | |\xi_1|^2 - |\xi_2|^2 |.$$
However, on $\Omega'_r$ we see from \eqref{alpha4} that
$$ -|\xi_1|^2  + |\xi_2|^2 - |\xi_3|^2 + |\xi_4|^2 = - 2 |\xi_{12}| |\xi_{14}| \cos \angle(\xi_{12},\xi_{14})
= O( N_3 N_1 \theta_0 )$$
and hence
$$ |\xi_1|^2 - |\xi_2|^2 = O( N_3 N_1 \theta_0 ) + O(N_3)^2$$
and the claim follows. 
\end{proof}

We now replace the resonance constraint $|\cos \angle(\xi_{12},\xi_{14})| \leq \theta_0$ with a simpler constraint.
Observe from elementary trigonometry that
$$ \angle(\xi_1, \xi_{14}) = O( \frac{N_4}{N_1} ); \quad \angle(\xi_3, \xi_{34}) = O( \frac{N_4}{N_3} )$$
and hence (since cosine is Lipschitz) 
$$ |\cos \angle(\xi_1, \xi_3)| = |\cos \angle(\xi_{14}, \xi_{34})| + O( \frac{N_4}{N_3} )
= |\cos\angle(\xi_{14}, \xi_{12})| + O( \frac{N_4}{N_3} ).$$
Thus on the resonance set $\Omega_r$ we have
$$ |\cos \angle(\xi_1, \xi_3)| \lesssim \theta_0 + \frac{N_4}{N_3}.$$
From these observations, we can now bound the left-hand side of \eqref{nmonster} by
$$ [m(N_1)^2 N_1 N_3 \theta_0 +  m(N_3)^2 N_3^2]
\| \tilde v^*_2 * \tilde v^*_4 * \tilde F \|_{L^\infty_{\tau,\xi}} $$
where
$$ F(t,x) := \int_\R \int_\R \int_{\R^2} \int_{\R^2} 1_{|\cos \angle(\xi_1, \xi_3)| \lesssim \theta_0 + \frac{N_4}{N_3}}
e^{i (t(\tau_1+\tau_3) + x \cdot (\xi_1+\xi_3))} \tilde v_1(\tau_1, \xi_1) \tilde v_3(\tau_3, \xi_3)\ d\xi_1 d\xi_3 d\tau_1 d\tau_3.$$
Applying Cauchy-Schwarz and \eqref{bil-est} we can thus estimate the left-hand side of \eqref{nmonster} by
$$
\lesssim [m(N_1)^2 N_1 N_3 \theta_0 +  m(N_3)^2 N_3^2] 
\frac{N_4^{1/2}}{N_2^{1/2}} \|v_2\|_{X^{0,b-}} \|v_4\|_{X^{0,b-}} \|F\|_{L^2_{t,x}}.
$$
Inserting this into \eqref{nmonster} and gathering terms using \eqref{n-relations-2}, we reduce to showing that
\begin{equation}\label{Ftarg}
\|F\|_{L^2_{t,x}} \lesssim
\frac{m(N_1)^2 N_1^{5/2-} m(N_3) N_4^{1/2} [N^{-2+} + \theta_0^{1/2} N^{-3/2+}]}
{m(N_1)^2 N_1 \theta_0 +  m(N_3)^2 N_3}
\|v_1\|_{X^{0,b-}} \|v_3\|_{X^{0,b-}}.
\end{equation}
Note that we have dropped $m(N_4)$ here because our assumptions imply
that $N_4 \lesssim N$ so $m(N_4) \thicksim 1$.
   
To proceed further we use

\begin{lemma}[Angularly refined bilinear Strichartz estimate]   
Let $0 < N_1 \leq N_2$ and $0 < \theta < \frac{1}{50}$.  Then for any $v_1,v_2\in X^{0,1/2+}$ with spatial frequencies $N_1,N_2$ respectively, the spacetime function
$$ F(t,x) := \int_{\R^2} \int_{\R^2}  
e^{i (t(\tau_1+\tau_2) + x \cdot (\xi_1+\xi_2))} 
1_{|\cos \angle(\xi_1,\xi_2)| \leq \theta} \tilde v_1(\tau_1, \xi_1) \tilde v_2(\tau_2, \xi_2)\ d\xi_1 d\xi_2$$
obeys the bound
\begin{eqnarray}\label{eq:fr-Bourgain}
\|F\|_{L^2_{t,x}}
\lesssim \theta^{1/2} \|v_1\|_{X^{0,1/2+}} \|v_2\|_{X^{0,1/2+}}.
\end{eqnarray}
\end{lemma}

\begin{proof} 
We may assume that $\theta \ll N_1/N_2$ since the claim follows from \eqref{bil-est} otherwise.
By standard averaging arguments (see e.g. \cite[Lemma 2.9]{tao}) it suffices to prove the claim for
$L^2$ free solutions, or more precisely that if $\phi_1, \phi_2 \in L^2_x$ have spatial frequency $N_1,N_2$ respectively and
$$ F(t,x) := \int_{\R^2} \int_{\R^2} 1_{|\cos \angle(\xi_1,\xi_2)| \leq \theta} e^{-it(|\xi_1|^2+\xi_2|^2)} e^{ix \cdot(\xi_1+\xi_2)} 
\hat \phi_1(\xi_1) \hat \phi_2(\xi_2)\ d\xi_1 d\xi_2$$
then
\begin{equation}\label{bourgain-1}
\|F\|_{L^2_{t,x}}
\lesssim  \theta^{1/2} \|\phi_1\|_{L^2_x}\|\phi_2\|_{L^2_x}.
\end{equation}

We first verify the estimate in the special case when the Fourier transform of $\phi_j$ is supported in
an angular sector $\{ \xi_j: \arg(\xi_j) = l_j \theta + O(\theta) \}$
of width $O(\theta)$ for $j=1,2$, where $l_1, l_2$ are arbitrary
integers $1 \leq l_1, l_2 \leq \frac{2 \pi}{\theta}$.  
Observe that the spacetime Fourier transform of $F$ is given by the formula
\begin{eqnarray*}
\tilde F(\tau,\xi) &=& \int_{\R^2} 1_{|\xi_1| \sim N_1} 1_{|\xi-\xi_1| \sim N_2}
1_{\arg(\xi_1) = l_1 \theta + O(\theta)} 1_{\arg(\xi-\xi_1) = l_2 \theta + O(\theta)}
1_{|\cos \angle(\xi_1,\xi-\xi_1)| \leq \theta}\\
& & \hat \phi_1(\xi_1) \hat \phi_2(\xi-\xi_1) \delta(|\xi_1|^2 + |\xi-\xi_1|^2 + \tau)\ d\xi_1.
\end{eqnarray*}
From the cosine rule
$$ |\xi|^2 = |\xi_1|^2 + |\xi-\xi_1|^2 + 2 |\xi_1| |\xi-\xi_1| \cos\angle(\xi_1,\xi-\xi_1)$$
and the hypothesis $N_1 \leq N_2$ we thus see that $\tilde F(\tau,\xi)$ is zero unless $|\xi| \sim N_2$, and
$0.9 |\xi|^2 \leq -\tau \leq 1.1 |\xi|^2$.  Thus
we may take absolute values followed by Cauchy-Schwarz and estimate
\begin{align*}
|\tilde F(\tau,\xi)|^2 &\leq 1_{|\xi| \sim N_2, 0.9 |\xi|^2 \leq -\tau \leq 1.1 |\xi|^2}
\left[\int_{\R^2} 1_{|\xi_1| \sim N_1} 1_{\arg(\xi_1) = l_1 \theta + O(\theta)} 1_{\arg(\xi-\xi_1) = l_2 \theta + O(\theta)}\right.\\
& \left.|\hat \phi_1(\xi_1)| |\hat \phi_2(\xi-\xi_1)| \delta(|\xi_1|^2 + |\xi-\xi_1|^2 + \tau)\ d\xi_1\right]^2 \\
&\leq \int_{\R^2} |\hat \phi_1(\xi_1)|^2 |\hat \phi_2(\xi-\xi_1)|^2 \delta(|\xi_1|^2 + |\xi-\xi_1|^2 + \tau)\ d\xi_1 \\
& \times \sup_{|\xi| \sim N_2, 0.9 |\xi|^2 \leq -\tau \leq 1.1 |\xi|^2}
\int_{\R^2} 1_{|\xi_1| \sim N_1} 1_{\arg(\xi_1) = l_1 \theta + O(\theta)} 1_{\arg(\xi-\xi_1) = l_2 \theta + O(\theta)}\\&
\delta(|\xi_1|^2 + |\xi-\xi_1|^2 + \tau)\ d\xi_1.
\end{align*}
Integrating this in $\tau$ and $\xi$ and using Plancherel's theorem, we see that to prove \eqref{bourgain-1} it
will suffice to show that
$$
\int_{\R^2} 1_{|\xi_1| \sim N_1} 
1_{\arg(\xi_1) = l_1 \theta + O(\theta)} 1_{\arg(\xi-\xi_1) = l_2 \theta + O(\theta)}
\delta(|\xi_1|^2 + |\xi-\xi_1|^2 + \tau)\ d\xi_1 \lesssim \theta
$$
or equivalently that
$$ | \{ \xi_1 \in \R^2: \arg(\xi_1) = l_1 \theta + O(\theta); 
\arg(\xi-\xi_1) = l_2 \theta + O(\theta); |\xi_1 - \xi/2| = r + O(\eps/N_2) \} |
\lesssim \theta \eps$$
whenever $|\xi| \sim N_2$ and $r \sim N_2$, and $\eps$ is sufficiently small.  But if $\xi_1$ is closer to $0$ than to $\xi$, the angular constraint $\arg(\xi-\xi_1) = l_2 \theta + O(\theta)$ restricts the circle $|\xi_1-\xi/2|=r$ to an arc of length $O(\theta)$; similarly if $\xi_1$ is closer to $\xi$ than to $0$ using the angular constraint
$\arg(\xi_1) = l_1 \theta + O(\theta)$.  The claim follows.

Now we establish the general case.  We can subdivide
$$ \phi_1 = \sum_{l_1} \phi_{1,l_1}; \quad \sum_{l_2} \phi_{1,l_2}$$
where $l_1, l_2$ range over the integers between $0$ and $2\pi/\theta$, and the Fourier transform of $\phi_{j,l_j}$
is supported in an angular sector $\{ \xi_j: \arg(\xi_j) = l_j \theta + O(\theta) \}$.  Then we have
\begin{align*}
 F(t,x)& = \sum_{l_1,l_2} 
 \int_{\R^2} \int_{\R^2} 1_{|\xi_1| \sim N_1} 1_{|\xi_2| \sim N_2} 1_{|\cos \angle(\xi_1,\xi_2)| \leq \theta}
 e^{ix\cdot (\xi_1+\xi_2)}
  e^{-it|\xi_1|^2}\\
  & \hat \phi_{1,l_1}(\xi_1) e^{-it|\xi_2|^2} \hat \phi_{2,l_2}(\xi_2)\ d\xi_1 d\xi_2.
  \end{align*}
But observe that if $|\cos \angle(\xi_1,\xi_2)| \leq \theta$ then $|\arg(\xi_1) - \arg(\xi_2)| = \pi/2 + O(\theta)$ 
or $|\arg(\xi_1) - \arg(\xi_2)| = 3\pi/2 + O(\theta)$
and hence $|l_2 - l_1| = \pi/(2\theta) + O(1)$ or $|l_2 - l_1| = 3\pi/(2\theta) + O(1)$.  From this, the triangle inequality, and the preceding computation we have
$$ \|F\|_{L^2_{t,x}}
\lesssim \sum_{l_1,l_2: |l_2-l_1| = \pi/(2\theta) + O(1) \hbox{ or } 3\pi/2\theta + O(1)}
\theta^{1/2} \|\phi_{1,l_1}\|_{L^2_x} \|\phi_{2,l_2}\|_{L^2_x}.$$
Observe that for each $l_1$ there are only $O(1)$ values of $l_2$ which contribute to this sum, 
and vice versa.  Thus by Schur's test
$$ \|F\|_{L^2_{t,x}}
\lesssim 
\theta^{1/2} (\sum_{l_1} \|\phi_{1,l_1}\|_{L^2_x}^2)^{1/2} (\sum_{l_2} \|\phi_{2,l_2}\|_{L^2_x}^2)^{1/2}$$
and the claim follows by Plancherel's theorem.
\end{proof}

\begin{remark}  In the regime $\theta \ll N_1/N_2$, the estimate is sharp, as can be seen by using (time-localized) free solutions of data $\phi_1, \phi_2$ whose Fourier transforms are indicator functions of the rectangles
$[N_1 - \theta N_2, N_1 + \theta N_2] \times [-\theta N_1, \theta N_1]$ and $[-\theta N_2, \theta N_2] \times [N_2 - \theta N_1, N_2 + \theta N_1]$ respectively; we omit the details.  Of course in the regime $\theta \gg N_1/N_2$ the estimate \eqref{bil-est} is superior.  It appears that similar estimates also hold if the angular constraint
$|\cos \angle(\xi_1,\xi_2)| \leq \theta$ is replaced with similar constraints such as $|\angle(\xi_1,\xi_2) - \alpha| \leq \theta$ for some $\alpha \gg \theta$, but we will not need such variants here.
\end{remark}

Applying this lemma, we reduce to showing that
\begin{equation}\label{man}
(\theta_0 + \frac{N_4}{N_3})^{1/2}
\lesssim 
\frac{m(N_1)^2 N_1^{5/2-} m(N_3) N_4^{1/2} [N^{-2+} + \theta_0^{1/2} N^{-3/2+}]}
{m(N_1)^2 N_1 \theta_0 +  m(N_3)^2 N_3}
\end{equation}

We establish \eqref{man} by splitting into some cases.

First, suppose that $N_3 \gtrsim \frac{N_4}{\theta_0}$, then the left-hand side of \eqref{man} can be
bounded by $O(\theta_0^{1/2})$.  If we bound the denominator on the right-hand side crudely from above by $O(m(N_1)^2 N_1)$, and discard the $N^{-2+}$ term in the numerator, we thus reduce to showing that
$$ \theta_0^{1/2} \lesssim \frac{ m(N_1)^2 N_1^{5/2-} m(N_3) N_4^{1/2} \theta_0^{1/2} N^{-3/2+} }{m(N_1)^2 N_1 }$$
which simplifies to
$$ 1 \lesssim N_1^{3/2-} m(N_3)N_4^{1/2} N^{-3/2+}.$$
But this is clear by estimating $m(N_3) \gtrsim m(N_1)$, $N_4^{1/2} \gtrsim 1$ and using $N_1 \gtrsim N$.

Henceforth we assume $N_3 \ll \frac{N_4}{\theta_0}$, so that the left-hand side of \eqref{man} is $N_4^{1/2} N_3^{-1/2}$, which allows us to rearrange \eqref{man} as
\begin{equation}\label{man2}
 m(N_1)^2 N_1 \theta_0 +  m(N_3)^2 N_3
\lesssim m(N_1)^2 N_1^{5/2-} m(N_3) N_3^{1/2} [N^{-2+} + \theta_0^{1/2} N^{-3/2+}].
\end{equation}
Let us first consider the bound for $m(N_1)^2 N_1 \theta_0$.  Discarding the $N^{-2+}$ factor on the right, we reduce to
$$ m(N_1)^2 N_1 \theta_0  \lesssim m(N_1)^2 N_1^{5/2-} m(N_3) N_3^{1/2} \theta_0^{1/2} N^{-3/2+}$$
which simplifies to
$$ \theta_0^{1/2} \lesssim N_1^{3/2-} m(N_3) N_3^{1/2} N^{-3/2+}.$$
But this is certainly true since $m(N_3) N_3^{1/2} \gtrsim 1$ and $\theta_0 \leq 1$.

Now we consider the bound for $m(N_3)^2 N_3$ in \eqref{man2}.  Discarding the $\theta_0^{1/2} N^{-3/2+}$ term, we reduce to
$$ m(N_3)^2 N_3 \lesssim m(N_1)^2 N_1^{5/2-} m(N_3) N_3^{1/2} N^{-2+}$$
which simplifies to
$$ m(N_3) N_3^{1/2} \lesssim m(N_1)^2 N_1^{5/2-} N^{-2+}.$$
Since $m(N_3) N_3^{1/2} \lesssim m(N_1) N_1^{1/2}$, we reduce to
$$ 1 \lesssim m(N_1) N_1^{2-} N^{-2+}$$
which is true since $N_1 \gtrsim N$.  This completes the proof of Proposition \ref{quatro}, and Theorem \ref{main} follows.




\end{document}